\documentclass[11pt, a4paper]{amsart}
\usepackage{amsbsy,amsfonts,amsmath,amssymb,amscd,amsthm, mathtools}
\usepackage{mathrsfs,dsfont}
\usepackage{graphicx}
\usepackage{mathtools}
\usepackage{xcolor}
\usepackage{lipsum}
\usepackage{todonotes}
\usepackage{empheq, fancybox}
\usepackage{wrapfig}

\newtheorem*{Main Theorem}{Main Theorem}
\newtheorem*{Main Tech Theorem}{Main Technical Theorem}
\newtheorem{Theorem}{Theorem}[section]
\newtheorem{Definition}[Theorem]{Definition}

\theoremstyle{definition}
\newtheorem{Remark}[Theorem]{Remark}

\newcommand{\bR}{{\mathbb R}}

\newcommand{\cI}{{\mathcal I}}

\DeclareMathOperator{\cl}{cl}

\begin{document}
\title[Uniform expansivity outside the critical neighborhood]
{Uniform expansivity outside the critical neighborhood in the quadratic family}
\date{March 23, 2015}

\author{Ali Golmakani}
\address{Ali Golmakani. Universidade Federal de Alagoas,
Campus A.C.\ Sim\~{o}es, Av.\ Lourival Melo Mota, s/n,
Cidade Universit\'{a}ria, Macei\'{o} Alagoas 57072-900, Brazil}
\email{aligolmakani@gmail.com}

\author{Stefano Luzzatto}
\address{Stefano Luzzatto. Abdus Salam
International Centre for Theoretical Physics (ICTP),
Strada Costiera 11, 34151 Trieste, Italy}
\email{luzzatto@ictp.it}

\author{Pawe\l{} Pilarczyk}
\address{Pawe\l{} Pilarczyk. Institute of Science and Technology Austria,
Am Campus 1, 3400 Klosterneuburg, Austria}
\email{pawel.pilarczyk@ist.ac.at}

\begin{abstract}
We use rigorous numerical techniques to compute a lower bound for the exponent of expansivity outside a neighborhood of the critical point for thousands of intervals of parameter values in the quadratic family. We compute a possibly small radius of the critical neighborhood, and a lower bound for the corresponding expansivity exponent outside this neighborhood, valid for all the parameters in each of the intervals. We illustrate and study the distribution of the radii and these exponents. The results of our computations are mathematically rigorous. The source code of the software and the results of the computations are made publicly available at \href{http://www.pawelpilarczyk.com/quadratic/}%
{http://www.pawelpilarczyk.com/quadratic/}.
\end{abstract}

\maketitle


\bigskip

We report on some rigorous numerical results concerning expansivity exponents for various parameter intervals  for the dynamics outside a critical neighborhood in the quadratic family of maps.  For completeness and in order to motivate our study, we start in Section 1 with some general facts on the quadratic family. The reader who is already familiar with the basic  background may wish to skip directly to Section 2. 


\section{Background and basic definitions}

Let \(  I  \) be a compact interval and \(  f\colon I \to I  \) a map. It is well known that the dynamics generated by the iterates of \(  f  \) can range from very simple to extremely complex. Two basic examples are the maps on the unit interval \(  I=[0,1]  \) given by \(  f(x)=x/10  \), for which the dynamics is very ``regular'' since all initial conditions converge under iteration to the fixed point at the origin, and \(  g(x)= 10 x  \) mod 1, which instead has an infinite number of distinct periodic orbits, an infinite number of dense orbits, and other hallmarks of what is generally described as ``chaotic'' or ``stochastic'' dynamics. It turns out (in a way that can be formalized) that these two quite distinct kinds of dynamical behavior are essentially due to the fact that in the first example the map is \emph{contracting}, i.e., the derivative is everywhere smaller than 1, and therefore nearby points get closer under iteration, while in the second example the map is \emph{expanding}, i.e., the derivative is everywhere greater than 1, and therefore nearby points are moving away under iteration. 

\begin{figure}[htbp]\label{fig:bif}
\setlength{\unitlength}{0.5\textwidth}
\begin{picture}(1.1,0.6)
\put(0.05,0.05){\includegraphics[width=1.0\unitlength,height=0.5\unitlength]{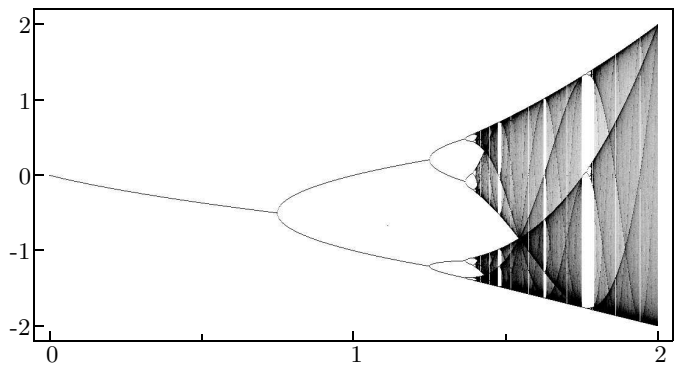}}
\linethickness{0.1mm}
\put(0.025,0.025){\line(1,0){1.05}}
\put(0.025,0.025){\line(0,1){0.55}}
\put(1.075,0.025){\line(0,1){0.55}}
\put(0.025,0.575){\line(1,0){1.05}}
\put(0.05,0.025){\line(0,1){0.015}}
\put(0.30,0.025){\line(0,1){0.01}}
\put(0.55,0.025){\line(0,1){0.015}}
\put(0.80,0.025){\line(0,1){0.01}}
\put(1.05,0.025){\line(0,1){0.015}}
\put(0.025,0.05){\line(1,0){0.015}}
\put(0.025,0.175){\line(1,0){0.015}}
\put(0.025,0.30){\line(1,0){0.015}}
\put(0.025,0.425){\line(1,0){0.015}}
\put(0.025,0.55){\line(1,0){0.015}}
\footnotesize
\put(0.045,-0.01){0}
\put(0.545,-0.01){1}
\put(1.045,-0.01){2}
\put(-0.015,0.035){-2}
\put(-0.015,0.160){-1}
\put(0.0,0.285){0}
\put(0.0,0.41){1}
\put(0.0,0.535){2}
\normalsize
\end{picture}
\setlength{\unitlength}{1pt}
\caption{Bifurcation diagram for the quadratic family.}
\end{figure}

Things get much more complicated in examples which combine some regions of contraction with some regions of expansion. One of the best known and most studied such examples is the real quadratic family of maps 
\[
f_{a}(x) = a- x^{2}.
\]
The fixed points of \(  f_{a}  \) are given by solutions of the equation \(  a-x^{2}=x  \) which gives \(  x=-1/2 \pm \sqrt{1+4a}/2  \). Thus \(  f_{a}  \) has two fixed points for every \(  a>-1/4  \). For \(  a>0  \) these lie on opposite sides of the origin and letting \(  p_{a}=-1/2-\sqrt{1+4a}/2  \) denote the negative fixed point, and letting \(  I_{a}=[p_{a}, -p_{a}]  \)  denote the compact interval defined by this fixed point, it is then easy to see that for every \(  a\in (0, 2]  \) the map 
\[  f_{a}\colon I_{a}\to I_{a}  \]
is well defined, i.e., \(  f_{a}  \) maps the compact interval \(  I_{a}  \) to  itself. The derivative of \(  f_{a}  \) is \(  f'_{a}(x)=-2x  \) and therefore it is easy to check that for positive values of \(  a  \) the derivative at the fixed point \(  p_{a}  \) satisfies \(  |f_{a}'(p_{a})|>1  \) whereas \(  f_{a}'(0)=0  \). Thus the map is neither fully contracting nor fully expanding. Nevertheless, it is easy to check that for small values of the parameter \(  a  \) the fixed point \(  q_{a}= -1/2+\sqrt{1+4a}/2  \), which lies in the interior of the interval \(  I_{a}  \), satisfies \(  |f_{a}'(q_{a})|<1  \) and is indeed an attractor for every point in the interior of \(  I_{a}  \). As \(  a  \) increases, however, the fixed point \(  q_{a}  \) loses its stability, the proportion of the interval \(  I_{a}  \) in which the map is expanding increases, and generally the dynamics gets significantly more complicated. Basic numerical simulations obtained by iterating a more or less arbitrary initial condition give rise to the well known \emph{bifurcation diagram}, see Figure \ref{fig:bif},
which was first observed in the 1970's and generally credited to Feigenbaum and which suggests certain regions of parameters for which the dynamics is basically still regular, with most initial conditions eventually attracted to some periodic orbit, interspersed by regions of parameter space in which the dynamics is very wild and chaotic.
Rigorous analytical results over the last couple of decades have gone a very long way in explaining and clarifying this bifurcation diagram, and revealed some initially unexpected and surprising facts. We restrict ourselves here to the parameter interval 
\[
\Omega=[1.4, 2]
\]
since this is the region in which the bifurcation diagram exhibits interesting phenomena. A first and remarkable observation which was already made as far back as the 1940's by Ulam and von Neumann \cite{UlaNeu47} is that the parameter \(  2  \) can be thought of as ``stochastic'' in a precise mathematical sense (more formally, the corresponding map \(  f_{a}  \) admits an ergodic invariant probability measure which is absolutely continuous with respect to Lebesgue measure). In 1981, Jakobson \cite{Jak81} and then in 1985 Benedicks and Carleson \cite{BenCar85}, showed that the set \(  \Omega^{+}\subset \Omega  \) of such stochastic parameters is  actually ``large'' in the sense that it has \emph{positive Lebesgue measure} (see also some generalizations in \cite{Ryc88, Tsu91, Rov93, PacRovVia98, Thu99, LuzTuc99, LuzVia00}). Interestingly, however, it is also ``small'' in the sense that it is  topologically \emph{nowhere dense}, a fact that follows from the  remarkable result that the set  \(  \Omega^{-}\subset \Omega  \) of ``regular'' parameters (for which almost every initial condition is eventually attracted to a periodic orbit) is \emph{open and dense} in \(  \Omega  \) \cite{GraSwi97, Lyu97a, Lyu97b} (see also generalizations in \cite{koz03, kozshest07}). Thus both regular and stochastic parameters have positive probability in \(  \Omega  \). It was proved in \cite{Lyu02} that in fact the union \(  \Omega^{+}\cup \Omega^{-}\subset \Omega  \) has \emph{full measure} in \(  \Omega  \) and these are therefore the only two dynamical phenomena that occur with positive probability, even though there are infinitely many parameters for which different, and sometimes extremely bizarre, dynamics does occur.

On a more quantitative level,  some attention has been devoted to the natural question of 
 how the measure is shared between the regular and stochastic parameters. Interestingly there seems to be no heuristic argument of any kind suggesting either that most parameters should be regular or that most parameters should be stochastic.
 Recent rigorous numerical estimates in \cite{TucWil09} confirm previous calculations in \cite{SimTat91} and show that for the related family \(  f_{a}(x)=ax(1-x)  \) and the corresponding relevant parameter interval \(  [2,4]  \) the set of regular parameters occupies at least \(  10\%  \) of all parameters. While the calculations have not been carried out for the quadratic family in the form that we consider here, it can be expected that they would yield similar results. Parameters in \( \Omega^+ \) cannot be computed directly since they belong to a Cantor set (which makes them theoretically un-decidable, see \cite{ArbMat04}), but it is possible to obtain explicit bounds for the measure of \( \Omega^+ \). Jakobson first set out the theoretical framework and posed the question of the possibility of obtaining explicit estimates in \cite{Jak01, Jak04}.  In \cite{LuzTak} it was shown that \(  |\Omega|\geq 10^{-5000}  \) and improvements on these estimates have been recently announced by Yu-Ru \cite{Hua11} and Shishikura \cite{Shi12}, but the bounds obtained are still very small. Some work in progress by the authors of the current paper aims to extend these estimates significantly.

One of the key conditions used in all proofs of the positivity of the Lebesgue measure of the set \( \Omega^+ \) of stochastic parameters is the so-called \emph{uniform expansivity outside the critical neighborhood} which will be defined precisely in the next section.  This property holds under certain general abstract conditions, but for the quantitative results mentioned above, which include numerical bounds on the measure of \( \Omega^+ \), it is necessary to have very explicit numerical estimates for the expansivity and the size of the critical neighborhood outside of which the expansivity holds. The main objective of this paper is to present the rigorous results of such computations for a relatively large number of small parameter intervals.


\section{Expansion outside the critical neighborhood}

As mentioned above, the dynamics of a map \(  f\colon I \to I  \) is very much influenced by whether \(  f  \) is expanding or contracting. Even if \( f \) is not expanding we can still formulate the notion of \( f \) being expanding in certain region of the phase space. 

\begin{Definition}
For \(  \lambda>0  \), \(  f\colon I \to I  \) is  \(  \lambda  \)-uniformly expanding outside  \(  \Delta\subset I \) 
if  there exists \( C   >0  \)  s.t.  for every \(  x\in I \) and every \(  n\geq 1  \)  with   \(  f^{i}(x)\notin\Delta  \) for all \(  i=0,1,\ldots,n-1  \), we have 
\begin{equation*}\label{eq:man}
|(f^{n})'(x)|\geq Ce^{\lambda n}.
\end{equation*}
\end{Definition}

Notice  that \(  f  \) is trivially \( \lambda \)-uniformly expanding outside \(  \Delta  \) if 
\(  |f'(x)|> e^\lambda >1  \) for all \( x\notin \Delta  \) (in which case one could also take \(  C=1  \)), and so in particular any expanding map is expanding outside any set \( \Delta \).  However, 
in many situations this property holds even if the derivative of \(  f  \) outside \(  \Delta  \) is very small in some regions.   
Indeed, a celebrated and important result of Ma\~n\'e \cite{Man85} says that if all periodic orbits of \(  f  \) are hyperbolic repelling (i.e., \(  |(f^{\ell})'(p)|>1  \) where \(  \ell  \) is the period of the periodic orbit \(  p  \)) then \emph{\(  f  \) is uniformly expanding outside \emph{any} neighborhood \(  \Delta  \) of the  critical points} (i.e., all the points \(  c  \) such that \(  f'(c)=0  \)), and even more generally, even if \(  f  \) has attracting periodic orbits then it is \emph{uniformly expanding outside any region \(  \Delta  \) as long as \(  \Delta  \) contains the critical points and the ``immediate basins'' of the attracting periodic orbits}.  This theorem has multiple consequences and corollaries and is arguably one of the most fundamental technical results in one-dimensional dynamics, on which many deep and important subsequent results are based. 

Ma\~n\'e's Theorem is an abstract ``existence'' result, and a natural question, and indeed the question we address and study in this paper, is that of an explicit numerical study of 
the size of the region \(  \Delta  \) and
the expansivity exponent \(  \lambda  \) in specific examples, in particular in the quadratic family, which we use as a case-study. The computational techniques required for a rigorous numerical investigation of this problem were developed in a previous paper \cite{DKLMOP08}. Here we further refine those techniques and apply them to a large-scale and systematic study of the entire relevant parameter space.  

In order to formulate our results, we introduce some definitions as follows. 
In the case of maps \(  f_{a}  \) belonging to the quadratic family, with \(  a\in \Omega  \), it is known that there can always be at most one attracting periodic orbit and that its immediate basin is a neighborhood of the critical point \(  c=0  \), and that all other periodic orbits are always hyperbolic repelling. This allows us to make the following definitions. For each \(  a\in \Omega  \), let 
\[
\delta_{a}:=\inf\{\delta>0: f_{a} \text{ is uniformly expanding outside } \Delta = (-\delta, \delta)\}.
\]
If \(  f_{a}  \) has a periodic attractor then \(  \delta_{a}>0  \), and either \(  \delta_{a}  \) or \(  -\delta_{a}  \) are boundary points of the immediate basin of attraction of the periodic orbit. In this case we can define 
\[
\lambda_{a}:=\sup\{\lambda:  f_{a} \text{ is \(  \lambda  \)-uniformly  expanding  outside }  \Delta_a = (-\delta_a,\delta_a)\}
\]
and, by the Theorem of Ma\~n\`e mentioned above, we always have \(  \lambda_{a}>0  \). 
If \(  f_{a}  \) does not have periodic attractors, then \(  \delta_{a}=0  \), 
i.e., \(  f_{a}  \) is uniformly expanding outside every neighborhood of the critical point. In this case we define 
\[
\lambda_{a}:=\lim_{\epsilon\to 0}\{\sup\{\lambda:  f_{a} \text{ is \(  \lambda  \)-uniformly  expanding  outside }  (-\epsilon, \epsilon)\}\}.
\]
In general, we obviously have \( \lambda_a\geq 0 \) and, surprisingly, for 
 almost every \(  a\in \Omega^{+}  \) we even have \(  \lambda_{a}>0  \) \cite{NowSan98}, meaning that \(  f_{a}  \) is in fact uniformly \(  \lambda  \)-expanding outside \emph{every} neighborhood of the critical point for some \(  \lambda  \) \emph{independent of the neighborhood}. Notice, however, that \(  f_{a}  \) is not \(  \lambda  \)-uniformly expanding for any \(  \lambda  \) if \(  \Delta  \) reduces just to the single critical point itself. 

We now extend these definitions to intervals of parameter values. 
For a closed non-trivial interval \(  \omega\subseteq\Omega  \), we let 
\[
\delta_{\omega}:=\sup_{a\in\omega}\{\delta_{a}\}=
\inf\{\delta>0: f_{a} \text{ is uniformly expanding outside } \Delta = (-\delta, \delta) \text{ for all }  a\in \omega\}.
\]
Since the set of regular parameters with attracting periodic orbits is open and dense, every non-trivial interval of parameters must contain some such parameter, and therefore we always have \(  \delta_{\omega}>0  \).
Given a radius $\delta$ of a critical neighborhood $\Delta = (-\delta, \delta)$, we let
\[
\lambda_{\omega} (\delta) := \inf\{\lambda>0:  f_{a} \text{ is \(  \lambda  \)-uniformly  expanding  outside }  \Delta = (-\delta, \delta) \text{ for all } a\in \omega\}.
\]

In the remaining part of the paper, we systematically analyze the parameter space \(  \Omega  \) for the quadratic family. We subdivide it into a large number of small subintervals, and obtain a collection of rigorous explicit upper bounds \( \bar\delta_\omega\geq  \delta_{\omega}  \) for \( \delta_\omega \) for each of the subintervals $\omega$, and lower bounds \(  \bar\lambda_\omega\leq \lambda_{\omega} (\bar\delta_{\omega})  \), which are expected to approximate \(  \lambda_{\omega} (\delta_{\omega})  \). Note, however, that \( \bar\lambda_\omega \) is, in general, neither an upper nor a lower bound for \( \lambda_{\omega} (\delta_{\omega}) \); this is due to the fact that the overestimate in calculating \( \bar\delta_\omega \) implies that a larger critical neighborhood is considered when determining $\bar\lambda_\omega$, which may increase the chance for computing a larger expansivity exponent than for $\delta_\omega$. Nevertheless, we \emph{prove} the following:

\medskip
\begin{center}
\emph{
For every \( a\in \omega \), the map $f_a$ is \( \bar\lambda_\omega \)-uniformly 
expanding outside \( (-\bar\delta_\omega, \bar\delta_\omega) \). 
}
\end{center}

\section{Numerical results and discussion}
\label{sec:results}

We now give plots which summarize the numerical results we obtained.  Details about the numerical algorithms and procedures will be given in Section~\ref{sec:procedures} together with various remarks about the computational issues involved.  We emphasize that all the results are absolutely rigorous. 

We subdivide the parameter interval \(  \Omega=[1.4, 2]  \) into $N = 60{,}000$ adjacent subintervals \(  \{\omega_i\}_{i = 0}^{N - 1}  \) of essentially equal size; see Section~\ref{sec:intervals} for the precise statement and the technical details of this subdivision. In order to simplify the notation, we omit the index $i$ when describing the computations done for each individual interval $\omega_i$.

\begin{figure}[htbp]
\setlength{\unitlength}{0.8\textwidth}
\begin{picture}(1.0,0.7)
\put(0,0.7){\includegraphics[height=\unitlength,angle=-90]{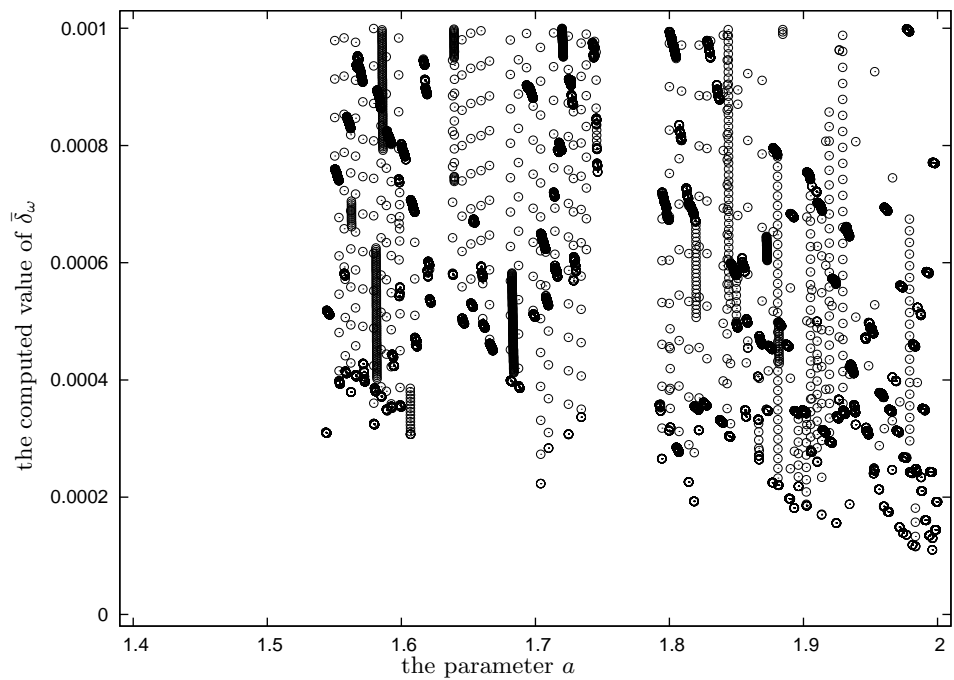}}
\footnotesize
\put(0.4,-0.01){the parameter $a$}
\put(0,0.2){\rotatebox{90}{the computed value of $\bar{\delta}_\omega$}}
\normalsize
\end{picture}
\setlength{\unitlength}{1pt}
\vskip 6pt
\caption{\label{fig:param2delta}
A possibly small radius $\bar{\delta}_\omega$ of the symmetric interval $\Delta_\omega$ around the origin for which it was possible to prove that $f_a$ is $\lambda$-uniformly expanding outside $\Delta_\omega$ with some $\lambda > 0$ for all $a \in \omega$. The circles indicate the pairs $(\omega, \bar{\delta}_\omega)$.}
\end{figure}

For each parameter interval \(  \omega  \), we (attempt to) compute a possibly tight upper bound $\bar{\delta}_\omega$ on $\delta_\omega$, as described in Section~\ref{sec:delta}. 
The computation was successful for $24{,}200$ out of the \( 60{,}000  \) parameter intervals yielding Figure \ref{fig:param2delta} 
which shows a plot of the value \(  \bar{\delta}_\omega  \) versus the parameter for all  intervals of parameters for which such computation was successful. Each of the $24{,}200$ data points is plotted as a dot surrounded by a small hollow circle. The circles tend to accumulate in some areas, and their superposition may occasionally give the impression of thick (black) dots. 

The lack of success of the calculation for a specific parameter interval \( \omega \) may be due to the fact that  \( \omega \) belongs to, or just intersects, some region of regular parameters with attracting periodic orbits with relatively large immediate basins, e.g., larger than the maximum size of \( 0.001 \) that we (somewhat arbitrarily) chose for the maximal allowed radius of the critical neighborhood. This is clearly the case, for instance, for the smaller parameters \( a\lesssim 1.55 \) and for the large obvious ``window'' around \( a\approx 1.77 \) corresponding to low period attracting periodic orbits which are  already obvious from the bifurcation diagram in Figure~\ref{fig:bif}. Although choosing a weaker constraint on the radius of the critical neighborhood (e.g., \( 0.1 \)) would considerably increase the number of successful calculations, in a longer perspective we are interested in further arguments regarding the measure of parameters for which ``stochastic'' dynamics occurs, and for this purpose it is absolutely necessary that the radius is very small. Moreover, expansivity in the phase space is much less interesting if it takes place very far from the critical point only. Therefore, by specifying the restrictive constraint on the radius of the critical neighborhood, we focus our attention on the results that are important from a wider perspective.

Parameter intervals for which the computation was successful, on the other hand, indicate a good deal of expansion, and therefore very probably a large proportion of stochastic parameters, especially if the expansion occurs outside some small critical neighborhoods. 
One perhaps noteworthy feature of Figure \ref{fig:param2delta}  consists of the several ``vertical'' rows of dots, occurring for example around \( a\approx 1.9 \). 
These seem to suggest the existence of a larger interval, formed by the union of a number of the smaller intervals we use in the computation, on all of which  we have expansion outside the largest critical neighborhood of radius \( 0.001 \), but such that some subintervals continue to show expansion even outside much smaller critical neighborhoods 
down to intervals of radius below \( 0.0002 \). Other curious patterns are the short diagonal patterns occurring for \( a\approx 1.7 \) for which we do not, however, have any heuristic explanation. 

\begin{figure}[htbp]
\setlength{\unitlength}{0.8\textwidth}
\begin{picture}(1.0,0.7)
\put(0,0.7){\includegraphics[height=\unitlength,angle=-90]{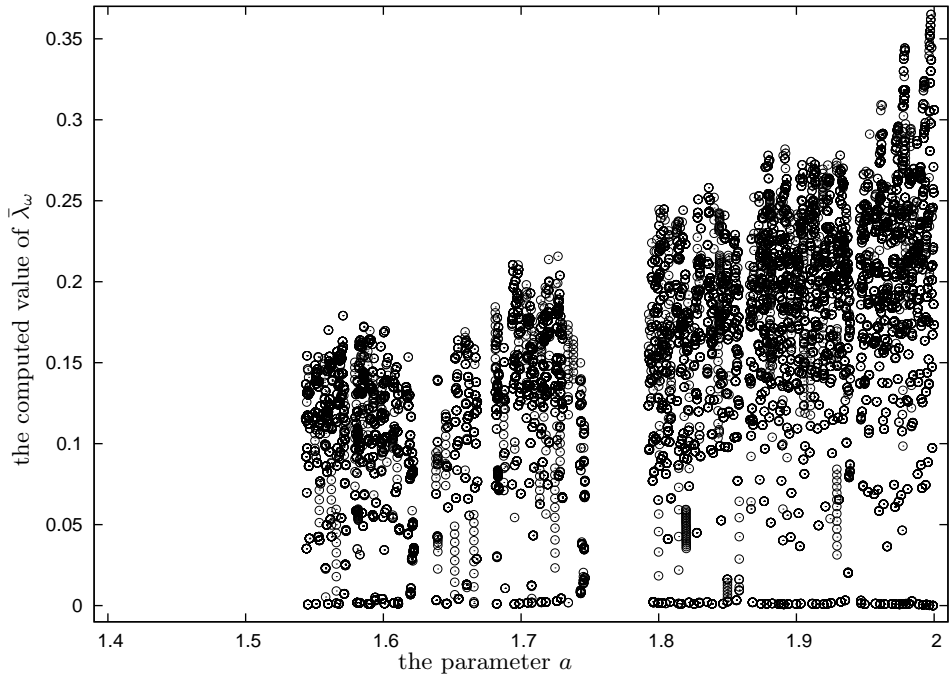}}
\footnotesize
\put(0.4,-0.01){the parameter $a$}
\put(0,0.2){\rotatebox{90}{the computed value of $\bar{\lambda}_\omega$}}
\normalsize
\end{picture}
\setlength{\unitlength}{1pt}
\vskip 6pt
\caption{\label{fig:param2lambda}
Values of $\bar{\lambda}_\omega$ for which it was possible to prove that $f_a$ is $\bar{\lambda}_\omega$-uniformly expanding outside $\bar{\Delta}_\omega = (-\bar{\delta}_\omega, \bar{\delta}_\omega)$ for all $a \in \omega$. The circles indicate the pairs $(\omega, \bar{\lambda}_\omega)$ for which the computation was successful.}
\end{figure}

For each parameter interval \(  \omega \) for which the computation of $\bar{\delta}_\omega$ was successful, we compute a lower bound $\bar{\lambda}_\omega$ for the expansion exponent $\lambda$ outside the critical neighborhood of radius $\bar{\delta}_\omega$, following the procedure described in Section~\ref{sec:graphs}, using a non-uniform subdivision of the complement of the critical neighborhood to $k = 20{,}000$ intervals (see Section~\ref{sec:k} for a more detailed description of the subdivision and a justification of this choice of the value of $k$). The results of this computation are summarized in Figure \ref{fig:param2lambda}.

Perhaps the main thing to notice here is how reasonably well distributed are the computed values of \( \bar\lambda_\omega \), without  any particular ``jumps'' or ``thresholds.'' We also notice that there are quite a few parameter intervals with particularly low values of \( \bar\lambda_\omega \), forming what looks almost like a horizontal line very close to 0. It seems that these are significant underestimates of the true values of \( \lambda_\omega \), as we discuss in more detail in Section~\ref{sec:k}. Indeed, we believe that most of these estimates represented in Figure \ref{fig:param2lambda} are underestimates due to the fact that improving the estimates requires a significant and major increase in required computational cost; we discuss this issue more in detail in the next section.


\begin{figure}[htbp]
\setlength{\unitlength}{0.8\textwidth}
\begin{picture}(1.0,0.7)
\put(0,0.7){\includegraphics[height=\unitlength,angle=-90]{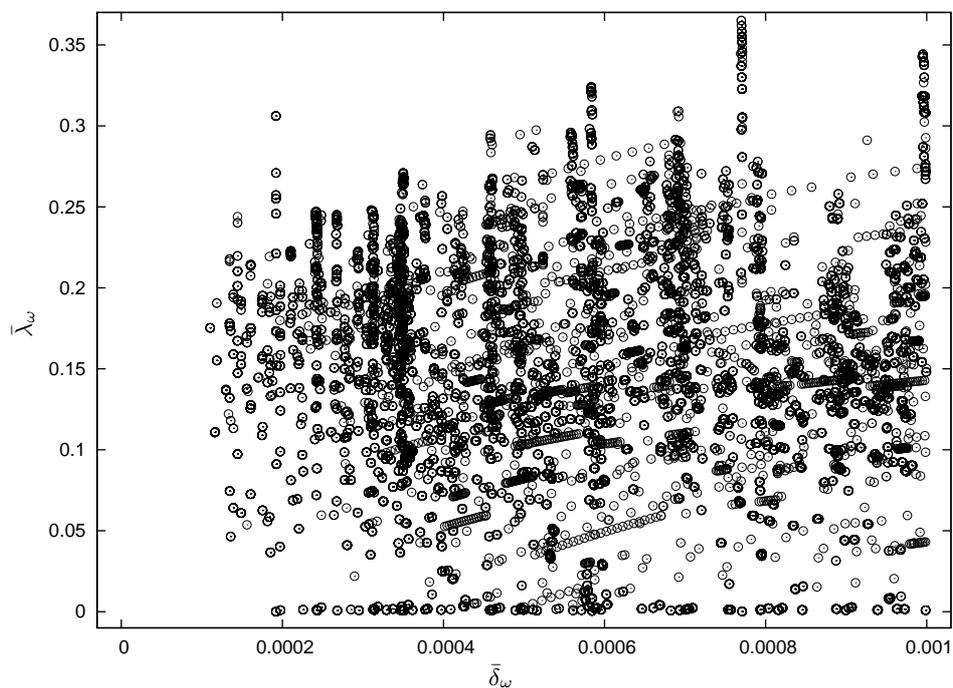}}
\footnotesize
\put(0.49,-0.02){$\bar{\delta}_\omega$}
\put(0,0.34){\rotatebox{90}{$\bar{\lambda}_\omega$}}
\normalsize
\end{picture}
\setlength{\unitlength}{1pt}
\vskip 6pt
\caption{\label{fig:delta2lambda}
The computed pairs of values $(\bar{\delta}_\omega, \bar{\lambda}_\omega)$ for those subintervals of $\Omega$ for which the computation was successful, indicated by circles.}
\end{figure}
%

An obvious question is whether there is any correlation between the computed values of \( \bar\lambda_\omega \) and those of \( \bar\delta_\omega \); a  naive conjecture might be that smaller critical neighborhoods imply a smaller expansivity exponent. This is, however, not the case, as shown in Figure~\ref{fig:delta2lambda} in which the values of \( \bar\lambda_\omega \) are plotted against the corresponding size \( \bar\delta_\omega \) of the critical neighborhood. Indeed, the picture is perhaps most remarkable in the total lack of correlation it shows between these two quantities, where both large and small values of \( \bar\lambda_\omega \) are essentially equally likely to correspond to relatively large or small critical neighborhoods. This is in fact consistent with the known qualitative theory of the quadratic family which, as mentioned above, says that for almost all stochastic parameters the expansivity constant is uniformly bounded from below for any \emph{arbitrarily small} critical neighborhood.

\begin{figure}[htbp]
\setlength{\unitlength}{0.9\textwidth}
\vskip -48pt
\begin{picture}(1.0,0.7)
\put(0,0.7){\includegraphics[height=\unitlength,angle=-90]{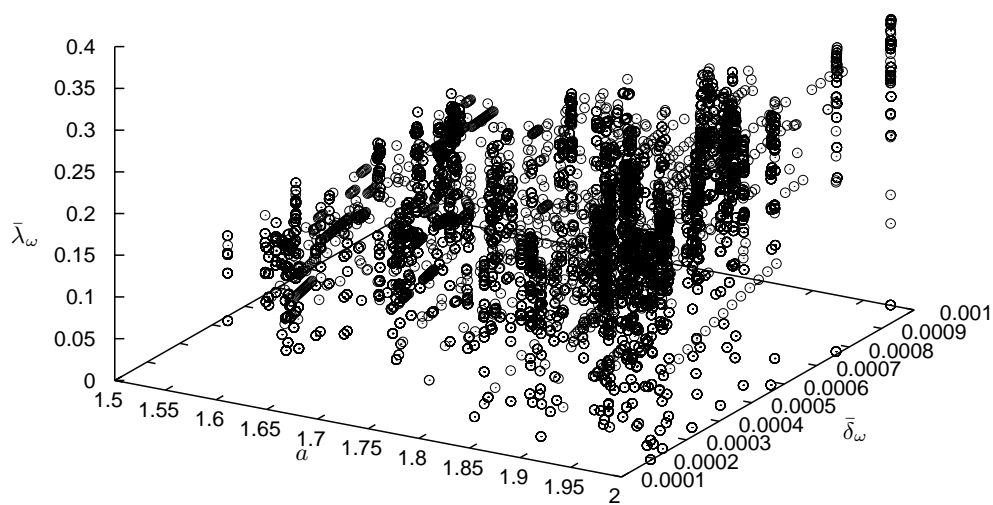}}
\footnotesize
\put(0.3,0.1){$a$}
\put(0.8,0.12){$\bar{\delta}_\omega$}
\put(0.04,0.3){$\bar{\lambda}_\omega$}
\normalsize
\end{picture}
\setlength{\unitlength}{1pt}
\vskip -24pt
\caption{\label{fig:3d3}
Three-dimensional plot of combined data of the two-dimensional plots in Figures \ref{fig:param2delta}, \ref{fig:param2lambda} and~\ref{fig:delta2lambda}.}
\end{figure}

We conclude this section with  a combined plot of all the three values: the parameter range and the corresponding values of \(  \delta  \) and \(  \lambda  \). This is shown in Figure \ref{fig:3d3}.

\section{Computational procedures}
\label{sec:procedures}

The computation of \(  \delta  \) and \(  \lambda  \) is essentially based on the procedures developed in \cite{DKLMOP08}, with some minor modifications, including in particular the use of a new and much more efficient algorithm in one part of the calculation, see discussion in Section \ref{sec:graphs} below. We clarify selected details and briefly outline the methods in this section and refer the interested reader to \cite{DKLMOP08} for more comprehensive explanations.

The source code of the specific software crafted for this paper
as well as raw results of the computations
are available at \href{http://www.pawelpilarczyk.com/quadratic/}%
{http://www.pawelpilarczyk.com/quadratic/}.
The computations can be run at multiple machines
(e.g., at a computer cluster) in a convenient way,
using the flexible dynamic parallelization scheme
introduced in~\cite{Pil10}, which is built into the software.


\subsection{Calculation of sub-intervals of $\Omega$}
\label{sec:intervals}

We use double-precision floating-point numbers following the ANSI/IEEE 754 standard, further called \emph{representable numbers}, which are available from the level of the C++ programming language, and implemented at the hardware level in modern processors.
For a real number $\alpha$, let $\widehat{\alpha}$ denote its representable approximation.
In order to subdivide the parameter interval \(  [\widehat{1.4}, 2]  \) into $N = 60{,}000$ adjacent subintervals \(  \{\omega_i\}_{i = 0}^{N-1}  \) of equal size (up to rounding), we first select the subdivision points. Namely, for each $i \in \{0, \ldots, N\}$, define
\[ \vartheta_i := \widehat{1.4} + \frac{(i / \gcd (i, N)) (2 - \widehat{1.4})}{N / \gcd (i, N)}, \]
where $\gcd (i, N)$ denotes the greatest common divisor of $i$ and $N$, and all the operations in the formula are done in double-precision floating-point arithmetic, with rounding the result of each operation to the nearest representable number (this is hardware-supported).
Then we set $\omega_i := [\vartheta_i, \vartheta_{i + 1}]$.
The purpose of truncating the fraction $i / N$ by the $\gcd (i, N)$ is that we would like to make sure
that the endpoints of the intervals computed for a finer subdivision, e.g., for $N = 120{,}000$,
agree with those for a coarser one, which need not be the case otherwise,
due to differences in rounding.

In order to avoid notation overload, in what follows, whenever explicit values are mentioned
for endpoints of an interval $\omega$, e.g., $[1.9999, 2]$,
we always mean their representable approximation.


\subsection{Computation of a lower bound for the expansion exponent $\lambda$}
\label{sec:graphs}

Given an interval $\omega \subset \Omega$ and a real number $\delta > 0$,
we compute a lower bound for $\lambda > 0$ such that
$f_a$ is uniformly expanding outside $\Delta := (-\delta, \delta)$ for all $a \in \omega$,
using a graph approach introduced in \cite{DKLMOP08}
with a non-uniform subdivision of the complement of the critical neighborhood
into a given number $k$ of intervals.
However, we introduce one major improvement.
Namely, we use a new algorithm \cite{Pil15}
for computing the minimum cycle mean in a directed graph,
instead of original Karp's algorithm \cite{Kar78}, 
which was applied in \cite{DKLMOP08}.
This tremendously reduces the memory usage, from $O (k^2)$ to $O (k)$.
For example, the previous algorithm
required some 400MB RAM for $k = 5{,}000$,
and with the new approach the memory usage drops down to some 12MB for the same graph.
This difference is more profound for large graphs: Over 50GB RAM required for $k = 50{,}000$
now reduces to some 20MB.
This improvement allows us to use much larger values of $k$
than it was possible previously in \cite{DKLMOP08},
and shifts the bottleneck from memory requirements to time constraints;
see the discussion in Section~\ref{sec:k}.

The basic idea of our approach to the computation
of a lower bound for $\lambda$ is to reduce the problem
to one of bounding the mean weights of paths in certain weighted digraphs
(directed graphs) related to the map \( f \) under consideration.
We let $G = (V,E,w)$ denote a \emph{weighted finite digraph},
where $V$ denotes the finite set of \emph{vertices},
$E \subset V \times V$ is the set of \emph{edges},
and $w \colon E \to \bR $ is the \emph{weight function}.
A {\em path}
is a nonempty finite sequence of edges
\(
\Gamma = (e_1, \ldots, e_n) 
\) 
such that \(
e_j= (v^0_j,v^1_j)\in E\ \text{and}\
v^1_j = v^0_{j + 1}.
\)
The path $\Gamma$ is called a \emph{cycle} if $v^0_1 = v^1_n$.
The \emph{weight} and \emph{mean weight}
of a path $\Gamma = (e_1, \ldots, e_n) $
are defined by
\[
W (\Gamma) = \sum_{j = 1}^{n} w (e_j)
\quad \text{ and } \quad
\overline{W} (\Gamma) = \frac{W(\Gamma)}{n},
\]
respectively.


In the setting of a quadratic map \( f_a \) as above, and a critical neighborhood \( \Delta \), we consider a collection of  intervals 
\( \cI = \{I_j \mid j = 1, \ldots, k\} \) with pairwise disjoint interiors and which cover \( I\setminus \Delta \), and say that the weighted digraph $G = (V, E, w)$ is a {\em representation of $f$ on $I \setminus \Delta$} provided that:
\begin{enumerate}
\renewcommand{\theenumi}{\alph{enumi}}
\item $V = \cI \cup \{\cl\Delta\}$;
\item\label{cond:edge}
\( \big\{e = (I_1, I_2) \in \cI \times V \mid
f (I_1) \cap I_2 \neq \emptyset \big\} \subset E \)
\item\label{cond:weight} For each $e = (I_1, I_2) \in E$,
\( 
w (e) \leq \inf \big\{ \log \big|Df (x)\big| :
x \in I_1 \cap f^{-1} (I_2) \big\}
\). 
\end{enumerate}

Observe the following straightforward
relationship between the weight of a path and the derivative along
points whose orbit is described by the path.
Given a point \( x \in I \setminus \Delta \) and a path
$\Gamma = (e_1, \ldots, e_n)$
such that $e_j = (I_{j-1},I_j)$ and
$f^{j} (x) \in I_j$ for all $j = 0, \dots, n$, we have
\begin{equation}\label{eqn:path}
\log \big|D f^n (x)\big| =
\sum_{j = 0}^{n - 1} \log \big|D f \big(f^j (x)\big)\big| \geq
W(\Gamma).
\end{equation}
The representation of the map and its derivatives to a weighted graph, and in particular the bound 
\eqref{eqn:path},
 reduce the problem of determining
expansion estimates to the computation of mean
weights of certain paths.
Specifically, the minimum mean weight of any cycle in $G$
provides a lower bound for the exponent $\lambda$ of interest.
This is the quantity computed by Karp's algorithm~\cite{Kar78}
or its improvement~\cite{Pil15}.

Some remarks were made in \cite{DKLMOP08} concerning the effect of using partitions with larger or smaller numbers of elements, or partitions with elements of variable size. As a result of those discussions, and a few additional tests, we choose a non-uniform partitioning of the interval $I_\omega := \bigcup_{a \in \omega} I_a$, as an apparently more effective one.

Note that the computations are done using rigorous numerics
(with controlled rounding directions), and are conducted
for the entire range $\omega$ of the parameters $a$ at a time
(by taking lower or upper bounds on the corresponding values, where appropriate).


\subsection{Computation of an upper bound on $\delta_\omega$}
\label{sec:delta}

Given an interval $\omega$, a possibly tight upper bound for $\delta_\omega$
is computed in the following way.
For $\delta > 0$, let $\bar{\lambda}_\omega^\delta$ denote the lower bound
for the expansion exponent outside $(-\delta, \delta)$,
computed with the procedure described in Section~\ref{sec:graphs}
for the parameter interval $\omega$,
using a coarse partition of $k = 1{,}000$ subintervals.
Starting with an apriori chosen $\delta_0 := 0.001$,
we first compute $\bar{\lambda}_\omega^{\delta_0}$.
If $\bar{\lambda}_\omega^{\delta_0} \leq 0$
then the computations for $\omega$ are considered to fail,
and no suitable upper bound on $\delta_\omega$ is reported.
Otherwise, a possibly small $\delta_1 \in (0, \delta_0]$ is found
for which $\bar{\lambda}_\omega^{\delta_1} > 0$,
using $20$ steps of the bisection method applied to $[0, \delta_0]$.
This number is further called $\bar{\delta}_\omega$
and constitutes a rigorous upper bound for $\delta_\omega$.

At this point we remark that in spite of the successful computation
of a positive expansivity exponent for $\bar{\delta}_\omega$
using the coarse partition of $k = 1{,}000$ subintervals,
in some rare cases it may happen that $\bar{\lambda}_\omega$
computed using a finer partition turns out to be negative.
This may be due to differences in the subdivision of the interval,
because of using a non-uniform partition, and should be treated
as a numerical artifact. If this happens, we also treat this situation
as a failure in the computation of $\bar{\delta}_\omega$
and $\bar{\lambda}_\omega$.


\subsection{The choice of the partition size $k$}
\label{sec:k}

As it was already discussed and shown in \cite{DKLMOP08},
the larger the number of partition elements $k > 0$ is chosen,
the higher the lower bound for the expansion coefficient $\lambda$ is found.
However, the time complexity of the algorithm for computing this bound
is $O(k^3)$, which practically means that the cost increases considerably
with the increase in $k$.

\begin{figure}[htbp]
\hfill
\setlength{\unitlength}{0.42\textwidth}
\begin{picture}(1.0,0.7)
\put(0,0.7){\includegraphics[height=\unitlength,angle=-90]{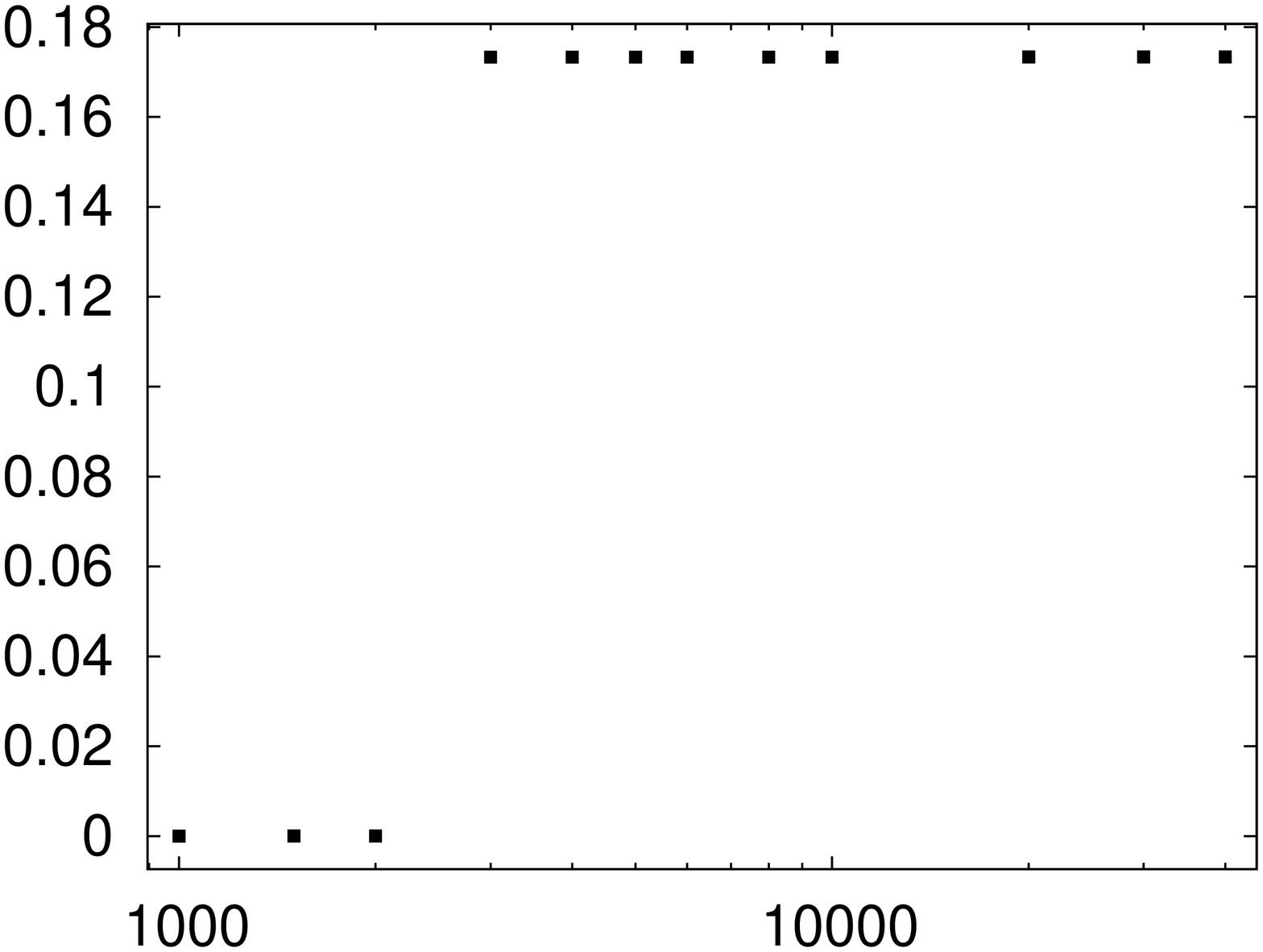}}
\footnotesize
\put(0.4,-0.02){partition size $k$}
\put(0,0.1){\rotatebox{90}{the computed exponent $\lambda$}}
\normalsize
\end{picture}
\setlength{\unitlength}{1pt}
\hfill
\setlength{\unitlength}{0.42\textwidth}
\begin{picture}(1.0,0.7)
\put(0,0.7){\includegraphics[height=\unitlength,angle=-90]{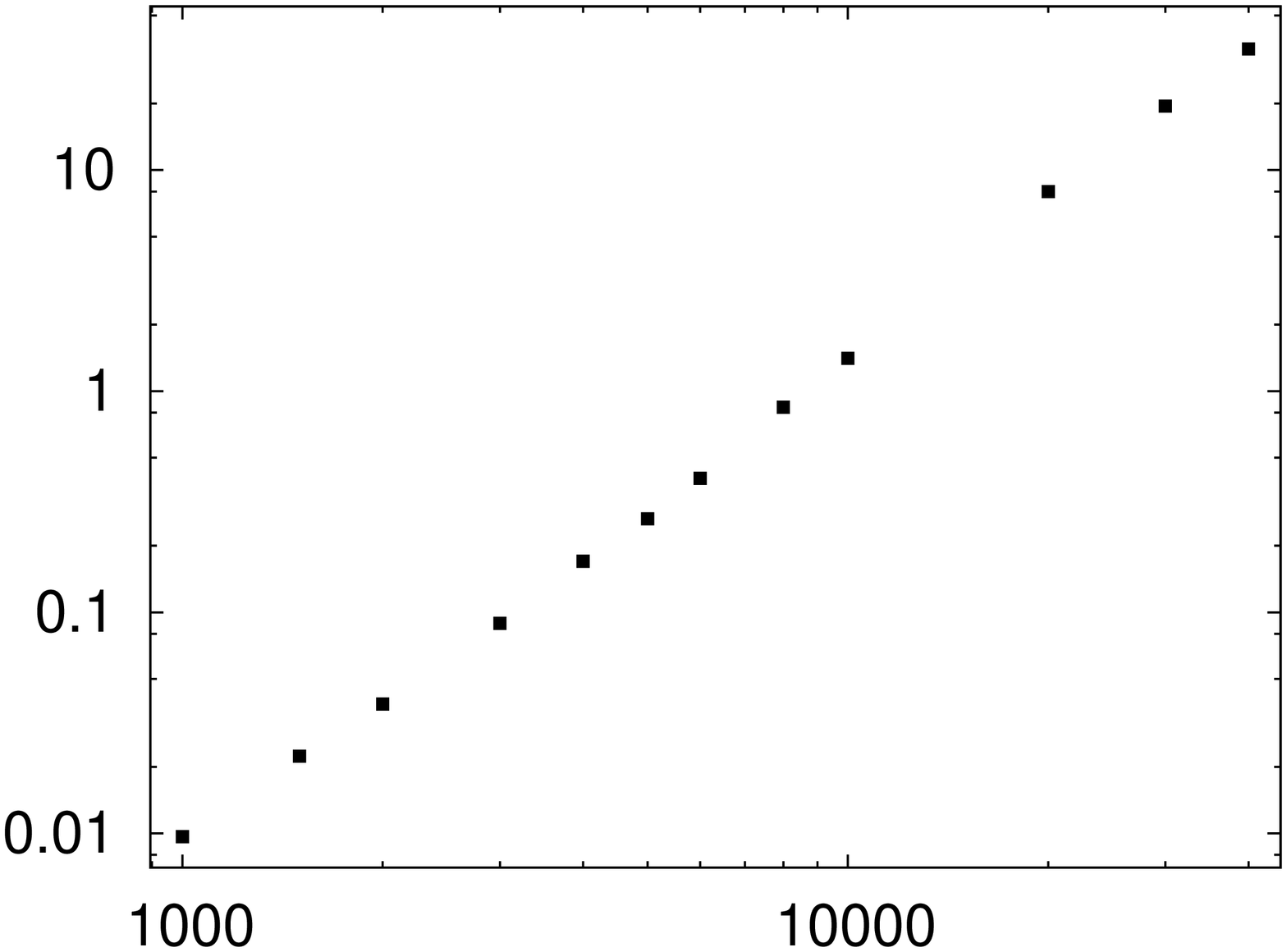}}
\footnotesize
\put(0.4,-0.02){partition size $k$}
\put(0,0.02){\rotatebox{90}{computation time $t$ [minutes]}}
\normalsize
\end{picture}
\setlength{\unitlength}{1pt}
\hfill
\vskip 6pt
\caption{\label{fig:part2lambdatime}
The exponent $\bar{\lambda}_\omega$ (the plot at the left-hand side)
and the time $t_k$ of computation of the exponent $\bar{\lambda}_\omega$
(the plot at the right-hand side)
as a function of partition size $k$,
computed for $\omega = [1.9999, 2]$.}
\end{figure}

In order to find a reasonable number $k$ for our comprehensive computations,
we conducted a test with $\omega = [1.9999, 2]$ and a selection
of different values $k$, ranging from $1{,}000$ to $40{,}000$.
The computed values of $\bar{\lambda}_\omega$
as a function of $k$ (note the logarithmic scale at the $k$ axis),
as well as the computation time (note the logarithmic scale at both axes)
are indicated in Figure~\ref{fig:part2lambdatime}.

\begin{figure}[htbp]
\setlength{\unitlength}{0.8\textwidth}
\begin{picture}(1.0,0.7)
\put(0,0.7){\includegraphics[height=\unitlength,angle=-90]{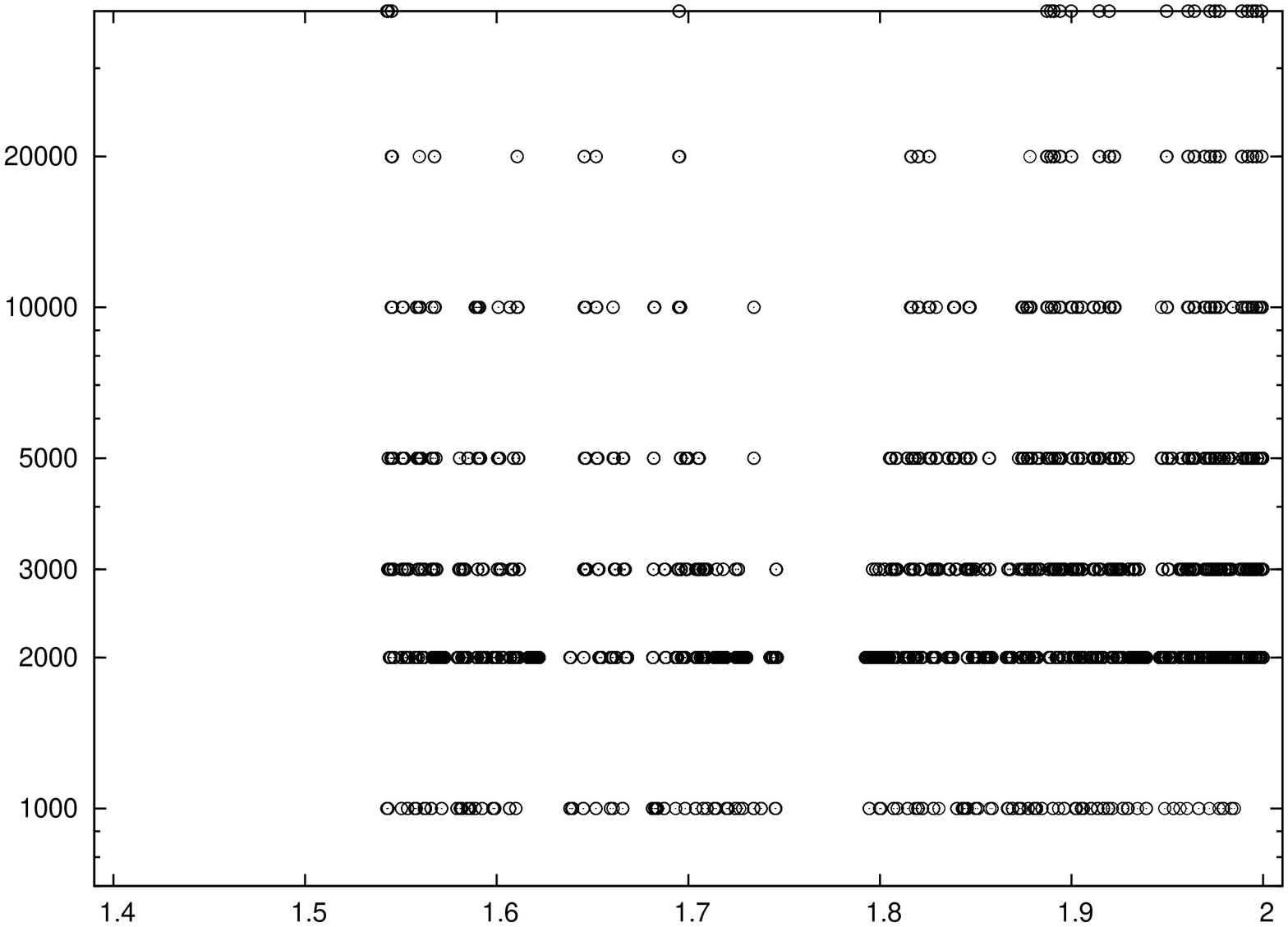}}
\footnotesize
\put(0.4,-0.01){the parameter $a$}
\put(0,0.2){\rotatebox{90}{the smallest $k$ for which $\bar{\lambda}_\omega > 0.001$}}
\put(0.98,0.092){(1{,}134)} 
\put(0.98,0.2){(13{,}196)} 
\put(0.98,0.263){(4{,}743)} 
\put(0.98,0.342){(2{,}921)} 
\put(0.98,0.455){(1{,}484)} 
\put(0.98,0.565){(387)} 
\put(0.98,0.67){(335)} 
\normalsize
\end{picture}
\setlength{\unitlength}{1pt}
\vskip 6pt
\caption{\label{fig:kLambdaThreshold}
The smallest value of $k \in \{1{,}000, 2{,}000, 5{,}000, 10{,}000, 20{,}000\}$
for which the computed $\bar{\lambda}_\omega$ is above $0.001$.
Circled dots at the top edge of the figure indicate the cases
in which the value of $\bar{\lambda}_\omega$ is below this threshold
for all the values of $k$ in the selection. The number of dots in each row
is indicated in parentheses at the right-hand-side of the diagram.}
\end{figure}

As one can see, the computed exponent $\bar{\lambda}_\omega$
is essentially $0$ up to a certain threshold,
and remains nearly constant afterwards.
It seems that taking too low a value for $k$ does not allow
to see the expansivity, which is then recovered after a certain threshold.
We conducted the computation of $\bar\lambda_{\omega}$
for $N = 60{,}000$ and a few values of $k$ between
$1{,}000$ and $20{,}000$; see Figure~\ref{fig:kLambdaThreshold}.
This experiment shows that the discussed threshold is different
for different parameter intervals.


Out of the $24{,}200$ subintervals for which the computation
of the expansivity exponent was successful, a value $\bar\lambda_{\omega} > 0.001$
is encountered in all but $335$ cases when $k$ reaches $20{,}000$,
which constitutes over $98.6\%$ of all the cases.
Taking this into consideration,
and also bearing in mind the cost of the computations
(which increases nonlinearly with the increase in $k$),
we decided to use $k = 20{,}000$ for all the plots discussed in Section~\ref{sec:results}.
Indeed, this value of $k$ seems to be beyond the threshold
that allows the computation of the seemingly correct exponent $\bar{\lambda}_\omega$
for a vast majority of intervals $\omega$,
yet the computation cost is reasonable
and takes about $4$ days on a relatively modern server
capable of running $32$ processes at its full speed.
Therefore, it seems that there is no point in increasing $k$ any further.
This justifies our choice of $k$.


\section*{Acknowledgments}

AG and PP were partially supported by Abdus Salam
International Centre for Theoretical Physics (ICTP).
Additionally, research conducted by PP has received funding
from Fundo Europeu de Desenvolvimento Regional (FEDER)
through COMPETE -- Programa Operacional Factores de Competitividade (POFC)
and from the Portuguese national funds
through Funda\c{c}\~{a}o para a Ci\^{e}ncia e a Tecnologia (FCT)
in the framework of the research project FCOMP-01-0124-FEDER-010645
(ref.\ FCT PTDC/MAT/098871/2008);
and from the People Programme (Marie Curie Actions)
of the European Union's Seventh Framework Programme (FP7/2007-2013)
under REA grant agreement no.\ 622033.

The authors gratefully acknowledge the Department of Mathematics of Kyoto University
for providing access to their server for conducting the computations described in the paper.



\begin{thebibliography}{99}

\bibitem{ArbMat04} 
A. Arbieto and C. Matheus. 2004.
\newblock Decidability of Chaos for Some Families of Dynamical Systems.
\newblock \emph{Foundations of Computational Mathematics}  4 (3): 269--275.


\bibitem{BenCar85}
M. Benedicks and L. Carleson. 1985.
\newblock On iterations of $1 - ax^2$ on $(-1, 1)$.
\newblock \emph{Ann. Math.} 122: 1--25.


\bibitem{DKLMOP08}
S. Day, H. Kokubu, S. Luzzatto, K. Mischaikow, H. Oka and P. Pilarczyk. 2008.
\newblock Quantitative hyperbolicity estimates in one-dimensional dynamics.
\newblock {\em Nonlinearity} 21: 1967--1987.

\bibitem{GraSwi97}
J. Graczyk and G. Swiatek. 1997.
\newblock Generic hyperbolicity in the logistic family.
\newblock {\em Ann. Math.} 146: 1--52.

\bibitem{Hua11} 
Y.-R. Huang. 2011.
\newblock Measure of parameters with acim nonadjacent to the Chebyshev value in the quadratic family.
\newblock PhD Thesis. University of Maryland.

\bibitem{Jak81}
M. V. Jakobson. 1981
\newblock Absolutely continuous invariant measures for one--parameter families of one--dimensional maps
\newblock {\em Commun. Math. Phys.} 81: 39-88.

\bibitem{Jak01}
M. V. Jakobson. 2001.
\newblock Piecewise smooth maps with absolutely continuous invariant measures and uniformly scaled Markov partitions.
\newblock {\em Proc. Symp. Pure Math.} 69: 825--881.
 
\bibitem{Jak04}
M. V. Jakobson. 2004
\newblock Parameter choice for families of maps with many critical points.
\newblock {\em Modern Dynamical Systems and Applications}
(Cambridge: Cambridge University Press).

\bibitem{Kar78}
R. M. Karp. 1978.
\newblock A characterization of the minimum cycle mean in a digraph.
\newblock {\em Discrete Math.}, 23(3): 309--311.

\bibitem{koz03}
O. Kozlovski. 2003.
\newblock Axiom A maps are dense in the space of unimodal maps in the $C^k$ topology.
\newblock {\em Ann. Math.} 157: 1--43.

\bibitem{kozshest07}
O. Kozlovski, W. Shen and S. van Strien. 2007.
\newblock Density of hyperbolicity in dimension one.
\newblock {\em Ann. Math.} 166: 145--182.

\bibitem{LuzTuc99}
S. Luzzatto and W. Tucker. 1999.
\newblock Non-uniformly expanding dynamics in maps with singularities and criticalities.
\newblock Institut Des Hautes Etudes Scientifiques. Publications Math\'{e}matiques (89): 179--226 (2000).

\bibitem{LuzVia00}
S. Luzzatto and M. Viana. 2000.
\newblock Positive Lyapunov exponents for Lorenz-like families with criticalities.
\newblock Ast\'{e}risque (261): xiii, 201--237.

\bibitem{LuzTak}
S. Luzzatto and H. Takahashi. 2006.
\newblock Computable conditions for the occurrence of non-uniform hyperbolicity in families of one-dimensional maps.
\newblock {\em Nonlinearity} 19: 1657--1695.

\bibitem{Lyu97a}
M. Lyubich. 1997.
\newblock Dynamics of quadratic polynomials. I.
\newblock {\em Acta Math.} 178: 185--247.

\bibitem{Lyu97b}
M. Lyubich. 1997.
\newblock Dynamics of quadratic polynomials. II.
\newblock {\em Acta Math.} 178: 247--297.

\bibitem{Lyu02}
M. Lyubich. 2002.
\newblock Almost every real quadratic map is either regular or stochastic.
\newblock {\em The Annals of Mathematics} 156 (1): 1--78.

\bibitem{Man85} 
R. Ma\~n\'e. 1985.
\newblock Hyperbolicity, sinks and measure in one-dimensional dynamics.
\newblock {\em Communications in Mathematical Physics} 100 (4): 495--524.

\bibitem{NowSan98} 
T. Nowicki and D. Sands. 1998.
\newblock Non-uniform hyperbolicity and universal bounds for $S$-unimodal maps.
\newblock {\em Inventiones Mathematicae} 132 (3): 633--680.

\bibitem{PacRovVia98}
M. J. Pacifico, A. Rovella and M. Viana. 1998.
\newblock Infinite-modal maps with global chaotic behavior.
\newblock {\em Annals of Mathematics}. Second Series 148 (2): 441--484

\bibitem{Pil10}
P. Pilarczyk. 2010.
\newblock Parallelization method for a continuous property.
\newblock {\em Foundations of Computational Mathematics} 10 (1): 93--114.

\bibitem{Pil15}
P. Pilarczyk.
\newblock A space-efficient algorithm for computing the minimum cycle mean in a directed graph.
\newblock Submitted.

\bibitem{Rov93}
A. Rovella. 1993.
\newblock The dynamics of perturbations of the contracting Lorenz attractor.
\newblock \emph{Boletim Da Sociedade Brasileira De Matem\'atica}. Nova S\'erie 24 (2): 233--259.

\bibitem{Ryc88}
M. R. Rychlik. 1988
\newblock Another proof of Jakobson's theorem and related results.
\newblock {\em Ergod. Theory Dyn. Syst.} 8: 93--109.
 
\bibitem{SimTat91} 
C. Sim\'o and J. C. Tatjer. 1991.
\newblock Windows of attraction of the logistic map.
\newblock In: \emph{European Conference on Iteration Theory} (Batschuns, 1989), 335--342.
World Sci. Publ., River Edge, NJ.

\bibitem{Shi12} 
M. Shishikura. 2012.
\newblock A Proof of Jakobson's Theorem via Yoccoz puzzles and the measure of stochastic parameters.
\newblock Conference slides.

\bibitem{Thu99}
H. Thunberg. 1999
\newblock Positive exponent in families with flat critical point.
\newblock {\em Ergod. Theory Dyn. Syst.} 19: 767--807.

\bibitem{Tsu91}
M. Tsujii. 1993
\newblock Positive Lyapunov exponents in families of one-dimensional dynamical systems
\newblock {\em Invent. Math.} 111: 113--137.

\bibitem{TucWil09}
W. Tucker and D. Wilczak. 2009.
\newblock A rigorous lower bound for the stability regions of the quadratic map.
\newblock {\em Physica D} 238 (18): 1923--1936.

\bibitem{UlaNeu47}
S. Ulam and J. von Neumann. 1947.
\newblock On combination of stochastic and deterministic processes.
\newblock {\em Bull AMS} 53 (11): 1120.

\end{thebibliography}
\end{document}